\documentclass[10pt,twoside]{article}
\usepackage{graphicx}
\usepackage{amsmath}
\usepackage{Latex-document}

\title{\bf Continuous Averaging in\vskip -2mm Dynamical Systems\vskip 6mm}

\author{D. Treschev\vspace*{-0.5cm}\thanks{Department of Mechanics and Mathematics,
Moscow State University, Vorob'evy Gory,
Moscow 119899, Russia. E-mail: dtresch@mech.math.msu.su}}
\date{\vspace{-8mm}}

\begin{document}

\maketitle

\thispagestyle{first} \setcounter{page}{383}

\begin{abstract}

\vskip 3mm

The method of continuous averaging can be regarded as a combination
of the Lie method, where a change of coordinates is constructed as
a shift along solutions of a differential equation and the Neishtadt
method, well-known in perturbation theory for ODE in the presence
of exponentially small effects. This method turns out to be
very effective in the analysis of one- and multi-frequency averaging,
exponentially small separatrix splitting and in the problem
of an inclusion of an analytic diffeomorphism into an analytic flow.
We discuss general features of the method as well
as the applications.

\vskip 4.5mm


\end{abstract}

\vskip 12mm

\section{The method}
\label{aver}\setzero \vskip-5mm \hspace{5mm}

There are several problems in the perturbation theory,
of real-analytic ordinary differential equations (ODE),
where standard methods do not lead to satisfactory results.
We mention as examples the problem of an inclusion of a
diffeomorphism into a flow in the analytic set up, and
the problem of quantitative description of exponentially small
effects in dynamical systems. In this cases one of possible
approaches is an application of the continuous averaging method.
The method appeared as an extension
of the Neishtadt averaging procedure \cite{nei1}.
We begin with the description of the method.

Let us transform the system
\begin{equation}
\label{eq_ini}
  \dot z = \widehat u(z),
\end{equation}
by using the change of variables
$z\mapsto Z(z,s)$.
Here $z$ is a point of the manifold $M$, $\widehat u$ is a smooth vector
field on $M$, $s$ is a non-negative parameter, and the change
is defined as a shift along solutions of the equation\footnote
 {Such method of constructing a change of variables is called
  the Lie method. The corresponding Hamiltonian version is
  called the Deprit-Hori method.}
\begin{equation}
\label{eq_aux}
  dZ/ds = f(Z, s),\qquad
  Z(z,0) = z,\quad
  0\le s\le S.
\end{equation}

Let the change $z\mapsto Z$ transform (\ref{eq_ini}) to the following system:
\begin{equation}
\label{eq_new}
  \dot Z = u (Z,s).
\end{equation}

Differentiating (\ref{eq_new}) with respect to $s$, we have:
$$
  \dot f(Z,s) = u_s(Z,s)
       + \partial_f u(Z,s)\quad
  \mbox{or }\;
  u_s = [u,f].
$$
Here $\partial_f$ is the differential operator on $M$, corresponding to
the vector field $f$, the subscript $s$ denotes the partial derivative,
and $[\cdot,\cdot]$ is the vector commutator:
$[u_1,u_2]=\partial_{u_1}u_2-\partial_{u_2}u_1$.
Putting $f=\xi u$, where $\xi$ is some fixed linear operator,
we obtain the Cauchy problem
\begin{equation}
\label{eq_aver}
  u_s = -[\xi u,u],\qquad u|_{s = 0} = \widehat u.
\end{equation}
We call the system (\ref{eq_aver})
averaging. The equation $f=\xi u$ is crucial for our method.
The vector field $f$ is usually constructed as a series in
the small parameter and not as a result of an application to
$u$ of an operator $\xi$, chosen in advance.

A nonautonomous analog of (\ref{eq_aver}) can
be easily constructed. If $\widehat u$ depends explicitly on $t$
then $f = \xi u$ also depends on $t$ and (\ref{eq_aver})
should be replaced by the system
\begin{equation}
\label{sys:a_nonaut}
         u_s = (\xi u)_t  -  [\xi u,u],\qquad
  u|_{s = 0} = \widehat u(z,t).
\end{equation}

Properties of the averaging system can be illustrated by the following
example. Consider the non-autonomous real-analytic system
\begin{equation}
\label{sys_ini}
  \dot z = \varepsilon\widehat u(z,t),\qquad
  z\in M.
\end{equation}
Here $\varepsilon$ is a small parameter, $\widehat u$ is
$2\pi$-periodic in $t$. Let us try to weaken the
dependence of $\widehat u$ on time by the change
$z\mapsto Z$ (\ref{eq_aux}) with $f = \xi u$.
We put\footnote
{Such an operator $\xi$ is called the Hilbert transform.}
\begin{equation}
\label{def:xi}
    \xi u(z,t,s)
  = \sum_{k\in {\bf Z}} i\sigma_k u^k (z,s) e^{ikt},\qquad
    \sigma_k = \mbox{sign\,} k ,
\end{equation}
where $u^k$ are Fourier coefficients in the expansion
$
  u(z,t,s) = \sum_{k\in {\bf Z}}  u^k (z,s) e^{ikt}.
$

Equation (\ref{sys:a_nonaut}) takes the form
\begin{eqnarray}
\label{ave1}
& u^k_s = - |k| u^k + i\varepsilon\sigma_k [ u^0, u^k ]
               - 2i\varepsilon \sum_{l+m=k, m<0<l} [u^l, u^m], & \\
\nonumber
& u^k|_{s = 0} = \widehat u^k,\qquad k\in{\bf Z} . &
\end{eqnarray}

To have an idea of properties of this system, we skip in (\ref{ave1})
the last term. The equations
$$
  u^k_s = - |k| u^k + i\varepsilon\,\sigma_k [u^0,u^k], \qquad
  u^k|_{s = 0} = \widehat u^k,\qquad k\in{\bf Z}
$$
can be solved explicitly:
\begin{equation}
\label{f_a_solut}
  u^k = e^{-|k|s} \widehat u^k \circ g^{i\varepsilon\sigma_k s},
\end{equation}
where $g^s$ is the time-$s$ shift $z|_{t=0}\mapsto z|_{t=s}$ along
solutions of the system $\dot z =  \widehat u^0(z)$.

The complex singularities of the functions $\widehat u^k\circ  g^\zeta$
of the complex variable $\zeta$ prevent an unbounded continuation
of the solutions (\ref{f_a_solut}) to all the set of positive
$s$. Nevertheless, the functions (\ref{f_a_solut})
can be made exponentially small in $\varepsilon$ since
$s$ can be chosen of order $\sim 1/\varepsilon$.

If $\varepsilon$ is not small, the operator (\ref{def:xi}) can be used to
smooth out the dependence of $u$ on $t$. Indeed, for arbitrarily
small $s>0$ the Fourier coefficients $u^k$ in (\ref{f_a_solut})
decrease exponentially fast in $k$ even if $\widehat u$
is just continuous in $t$.

Certainly, these hewristic arguments cannot be regarded as a proof of
the fact that systems of type (\ref{ave1})
can be used for averaging or smoothing. Rigorous statements
and estimates are based on the majorant method. In this paper
we do not go into technical details, but present general ideas
and applications.

If the vector field $\widehat u$ belongs to some subalgebra $\chi$
in the Lie algebra of vector fields on $M$,  it is natural
to look for the change $z\mapsto Z$ from the corresponding
Lie group of diffeomorphisms. This means that
$f$ in (\ref{eq_aux}) should be taken from $\chi$.
The same remains reasonable in the non-autonomous case.
Note that the operator (\ref{def:xi}) is such that if
$u(z,t,s)\in\chi$ for any $t$, the vector field $\xi u$
also belongs to $\chi$ for any $t$.

\section{Applications}
\setzero\vskip-5mm \hspace{5mm }

\subsection{Fast phase averaging: one-frequency case}
\vskip-5mm \hspace{5mm}

Mathematical models of various physical processes use
systems of ODE which contain an angular variable
changing much faster than other variables in the system.
Taking the fast phase as a new time, we can rewrite
the equations in the form
\begin{equation}
\label{onefr}
  \dot z = \varepsilon\widehat u(z,t,\varepsilon),\qquad z\in M,
\end{equation}
where $M$ is the $m$-dimensional phase space of the system,
and $\varepsilon$ is a small parameter.
The vector field $\widehat u$ is assumed to be
smooth and to depend on time $2\pi$-periodically.

It is well known that by a change of the variables it is possible
to weaken the dependence of the system~(\ref{onefr})
on time. In particular, by using the standard averaging method,
it is easy to construct
a $2\pi$-periodic in $t$ change of the variables $z\mapsto z_*$
such that the equations (\ref{onefr}) take the form
\begin{equation}
\label{averaged}
  \dot z_* = \varepsilon\widehat u^0(z_*)
              + \varepsilon^2\widehat u_*(z_*,\varepsilon)
                     + \varepsilon\tilde u(z_*,t,\varepsilon).
\end{equation}
Here the only term in the right-hand side depending explicitly
on time is $\varepsilon\tilde u = O(\varepsilon^K)$.
The natural $K$ is arbitrary and
$
 \widehat u^0(z) = \frac 1{2\pi} \int_0^{2\pi} \widehat u(z,t,0)\,dt.
$

Now suppose that $\widehat u$ is real-analytic in $z$.
Poincar\'e noted in some example that in this case power
series in $\varepsilon$ presenting a change of variables
eliminating time from the equations, exist but diverge:
terms at $\varepsilon^k$ in these series have the order $k!$.
In a general situation this statement has been proved by
Sauzin \cite{Sauz_Gevr}.

Neishtadt \cite{nei1} noted that in this case it is possible to obtain
in (\ref{averaged})
\begin{equation}
\label{exp_sm_pert}
  \tilde u = O(e^{-\alpha / \varepsilon}),\quad
  \alpha = \mbox{const} > 0
\end{equation}
($\varepsilon$ is assumed to be nonnegative).
The method Neishtadt used to prove this assertion
is based on a large (of order $1/\varepsilon$) number of
successive changes of variables. These changes weaken
gradually explicit dependence of the equations on time.
Ramis and Schafke \cite{Ram_Scha} obtained analogous results
analyzing diverging series, produced by the standard averaging method.

It is known also that in general a constant $A > \alpha$ exists
such that it is impossible to construct $2\pi$-periodic in $t$
change $z\mapsto z_*$, such that $\tilde u = O(e^{-A / \varepsilon})$.
This statement follows, for example,
from an estimate of the separatrix splitting rate in Hamiltonian
systems of type (\ref{onefr}) with one and a half degrees of freedom.

In this section we
estimate a ``maximal'' $\alpha$ for which the estimate
(\ref{exp_sm_pert}) is possible.

Suppose that the manifold $M$ is real-analytic. We fix its complex
neighborhood $M_{\bf C}$ and denote by $g^t$ the phase flow of the
averaged system\footnote
 {It would be more correct to call it by the first approximation
 averaged system, written with respect to the fast time.}
$\dot z = \widehat u^0(z)$.

Let $Q$ be a compact in $M$ and $V_Q$ its neighborhood in $M_{\bf C}$.
Suppose that for any real $s$ such that $|s| < \alpha$ and
for any point $z\in V_Q$ the map $g^{is}$ is analytic at
$z$ and moreover, $g^{is}(z)\in M_{\bf C}$. We define the set
$$
  U_{Q,\alpha} = \bigcup_{- \alpha < s < \alpha} g^{is} (V_Q).
$$

{\bf Theorem 1} \cite{Tre_RCD}. {\it Let the positive constants
$\alpha,\rho,\varepsilon_0$ be such that

(1) $U_{Q,\alpha}\subset M_{\bf C}$.

(2) The vector field $\widehat u$ is analytic in $z$ and $C^2$-smooth
    in $t,\varepsilon$ on
    $U_{Q,\alpha} \times {\bf T} \times [0,\varepsilon_0]$.

Then for sufficiently small $\varepsilon_0$, there exists a
$2\pi$-periodic in $t$ real-analytic in $z$ map
$ F: V'_Q \times {\bf T} \times (0,\varepsilon_0) \to M_{\bf C},\quad
  Q\subset V'_Q\subset V_Q,
$
such that

(a) The set $V'_Q$ is open in $M_{\bf C}$,

(b) $F(z,t,\varepsilon) = z_* = z + O(\varepsilon)$,

(c) $F$ transforms (\ref{onefr}) into (\ref{averaged})
   and the following estimate holds:
\begin{equation}
\label{tildev_global}
  |\tilde u(z,t,\varepsilon)| \le C e^{-\alpha / \varepsilon},\qquad
  z\in V'_Q,\quad
  t\in \Sigma_\rho,\quad
  \varepsilon\in [0,\varepsilon_0).
\end{equation}
}

Theorem~1 means in particular, that in the case when
components of the field $\widehat u$ are entire functions of
$z$, the quantity $\alpha$ in (\ref{tildev_global})
can be arbitrary positive number such that for all
$s\in [-\alpha,\alpha]$ the maps $z\mapsto g^{is}(z)$
are holomorphic at any point $z\in Q$.

Proof of Theorem~1 is based on the continuous averaging.
Namely, we solve the Cauchy problem~(\ref{sys:a_nonaut}),
with $\xi$ defined by (\ref{def:xi}).
The required change of variables corresponds to the value
$s = \alpha / \varepsilon$. The averaging can be performed inside a
subalgebra $\chi$ in the Lie algebra of all real-analytic vector fields
on $M$. In particular, if the initial vector field $\widehat u$
is Hamiltonian then $u_* = u(z,t,\varepsilon,\alpha / \varepsilon)$
is also Hamiltonian, and $F$ is symplectic.

\subsection{Averaging: multi-frequency case}
\vskip-5mm \hspace{5mm}

Consider a real-analytic slow-fast system
\begin{equation}
\label{sys2}
 \dot x = \omega + \varepsilon (\overline u(y)
                 + \widehat u(x,y,\varepsilon)),\quad
 \dot y =           \varepsilon (\overline v(y)
                 + \widehat v(x,y,\varepsilon)),\qquad
 x\in{\bf T}^n,\;  y\in{\bf R}^m.
\end{equation}
Average in $x$ of $\widehat u$ and $\widehat v$ is assumed to be
$O(\varepsilon)$. The frequency vector
$\omega\in{\bf R}^n$ is constant and non-resonant.\footnote
 {The case of non-constant frequencies ($\omega=\omega(y)$) can
  be reduced to this one in a small neighbourhood of an unperturbed
  (may be, resonant) torus $\{y = y^0\}$.
 }

We try to weaken the dependence of the right-hand side of the
system on the fast variables $x$ by a near-identity change
$(x\bmod 2\pi, y)\mapsto (x_\bullet \bmod 2\pi, y_\bullet)$.
We put
$$
              z = \Big(\begin{array}{c} x \\ y\end{array} \Big),\quad
\overline\omega = \Big(\begin{array}{c}\omega \\ 0\end{array} \Big),\quad
    \overline w = \Big(\begin{array}{c}\overline u\\
                                        \overline v \end{array}\Big),\quad
         \widehat w = \Big(\begin{array}{c}\widehat u\\ \widehat v\end{array}\Big).
$$
Note that the vector fields $\overline\omega$ and $\overline w$ commute.
The system (\ref{sys2}) takes the form
\begin{equation}
\label{sys3}
  \dot z = \overline\omega + \varepsilon (\overline w + \widehat w).
\end{equation}

{\bf A1. Diophantine condition.} { We assume that $n\ge 2$ and
the frequency vector $\omega$ is Diophantine}: there exist
$\gamma_0,\gamma > 0$ such that
\begin{equation}
\label{ineq:dioph}
  |\langle k,\omega\rangle| \ge \gamma_0 ||k||^{-\gamma}, \quad
  \mbox{for any }\,
  k \in {\bf Z}^n\setminus\{0\}.
\end{equation}

To formulate the next assumption, we need some definition.
Let $g^t$ be the phase flow of the system
\begin{equation}
\label{sys:ave}
  \dot z = \overline w(y).
\end{equation}

For any real-analytic function $f(y)$
with values in ${\bf C}^{n+m}$ and the vector $k\in{\bf Z}^n$ we put
$f_k = f(y) e^{i\langle k,x\rangle}$. The function
$$
  {\bf g}^s_k f= e^{-i\langle k,x\rangle} g_*^{-is} (f_k\circ g^{is}),
  \qquad s\in{\bf C}
$$
does not depend on $x$. Here $g_*^t$ is the differential of the map $g^t$.
Note that ${\bf g}^s_k$, $s\in{\bf R}$ include shifts along $g^t$ in purely
imaginary direction.
We put
\begin{eqnarray*}
 \Sigma_q &=& \{ x\in{\bf C}^n/(2\pi{\bf Z})^n : |\mbox{Im}\, x_j | \le q,\quad
                                            j = 1,\ldots,n \}, \\
      V_\nu    &=& \{ y\in{\bf C}^m :
                \mbox{Re}\, y \in (\overline{\cal D}+\nu)\subset{\bf R}^m, \quad
                                          |\mbox{Im}\, y_l | \le \nu, \quad
                                            l = 1,\ldots,m \},
\end{eqnarray*}
where $\overline{\cal D}$ is a compact domain
and $\overline{\cal D} + \nu$ is the $\nu$-neighborhood of
$\overline{\cal D}$.

Expand the function $\widehat w$ into the Fourier series:
$$
    \widehat w(z,\varepsilon)
  = \sum_{k\in{\bf Z}^n} \widehat w^k(y,\varepsilon) e^{i\langle k,x\rangle}.
$$

{\bf A2. Analyticity.} { Let the constant $\alpha > 0$ be such
that for any real $s \in [-\alpha,\alpha]$, for any $k\in{\bf Z}$
and for any $(x,y) \in \overline\Sigma_q\times\overline V_\nu
        = \mbox{\rm closure}(\Sigma_q)\times\mbox{\rm closure}(V_\nu)$
the map $(x,y)\mapsto g^{is}(x,y)$ is analytic in $x,y$ and
the map $y\mapsto {\bf g}_k^{s}\widehat w^k(y)$ is analytic in $y$.

Moreover, suppose that for some real $\rho$
$$
      |{\bf g}_k^s
      \widehat w^k |_{\overline V_\nu \times [0,\varepsilon_0)}
  \le \mu ||k||^{-\rho} e^{-q ||k||},\qquad
      k\ne 0.
$$

The function $\widehat w^0$ is of order $\varepsilon$ and we assume that}
$
    |{\bf g}_0^s \widehat w^0 |_{\overline V_\nu \times [0,\varepsilon_0)}
 \le\varepsilon \mu_0 .
$

The constants $\mu,\mu_0,\nu,q,\rho,\gamma$ must satisfy
some conditions \cite{PT2000}. Here we replace these conditions by more
restrictive, but simple ones.

{\bf A3}. $\mu,\mu_0,\nu,q,\rho,\gamma$ are positive, do not
depend on $\varepsilon$, and $\rho > \gamma/(\gamma+1)$.

{\bf Theorem 2} \cite{PT2000}. {\it Suppose that assumptions {\bf
A1}--{\bf A3} hold. Then there exists a change of variables
\begin{equation}
\label{change:final}
  z \mapsto z_\bullet = f(z,\varepsilon),\quad
  f: \Sigma_{2q/3} \times V_{2\nu/3} \times [0,\varepsilon_0)
    \to  {\bf C}^n / (2\pi{\bf Z}^n) \times {\bf C}^m
\end{equation}
such that
$f$ is analytic in $z$, smooth in $\varepsilon$,
$f(z,\varepsilon) = z + O(\varepsilon)$, and
the system {\rm (\ref{sys3})} takes the form
  \begin{equation}
  \label{sys:final}
    \dot z_\bullet = \overline\omega + \varepsilon (\overline w(y_\bullet)
                                       + w_\bullet(z_\bullet,\varepsilon)).
  \end{equation}

Let $w^0_\bullet(y_\bullet,\varepsilon)$ be the average in $x_\bullet$
of $w_\bullet$. Then
$w^0_\bullet(y_\bullet,0) = 0$.
Moreover,
\begin{equation}
\label{est:main1}
      |w_\bullet - w^0_\bullet|
  \le C\mu \varepsilon^{\rho / (\gamma + 1)}
        e^{-\overline q \varepsilon^{-1/(\gamma + 1)}} ,\qquad
      z_\bullet\in \Sigma_{q/2} \times V_{\nu/2},
\end{equation}
where $C$ is a constant, not depending on $\varepsilon$ and $\mu$,
\begin{equation}
\label{barq=}
  \overline q = (1 + \gamma^{-1})
                    (\gamma\gamma_0\alpha q^\gamma)^{1/(\gamma + 1)}.
\end{equation}
}

If the system {\rm (\ref{sys3})} is Hamiltonian with respect to a certain
symplectic structure $\Omega$ then {\rm (\ref{sys:final})}
is also $\Omega$-Hamiltonian and the change
{\rm (\ref{change:final})} is $\Omega$-symplectic.

In \cite{PT2000} we also present another theorem which shows
that the Fourier series $w_\bullet - w_\bullet^0$ can be divided
into 2 parts: one is small and for another we have a sort
of control.

Results analogous to Theorem 2 (without estimates for $\alpha$)
are contained in \cite{Simo,Bambu,LMS}.

\subsection{Exponentially small separatrix splitting}
\vskip-5mm \hspace{5mm}

The phenomenon of exponentially small separatrix splitting was discovered
by Poincar\'e \cite{poi}. Intensive quantitative studying
of the problem was initiated by papers \cite{hms, laz_0}
(see also \cite{delsea}). The method proposed by Lazutkin
with collaborators \cite{laz_1, laz_2, Vasya, Vasya2}
is based on an analysis of the separatrices in the complex domain.
Another method was used in \cite{G3}, where direct
expansions of the Poincar\'e-Melnikov integral in an additional
parameter are analyzed.
The resurgent analysis is applyed to these problems in \cite{Sauzin}.

The main difficulty of the problem is that the traditional
Poincar\'e-Melnikov method can not be applied directly.
Indeed, its error has the order of square of the perturbation.
Hence, the error considerably exceeds the expected result.

Exponentially small separatrix splitting
can be studied with the help of the continuous averaging method.
The main idea of to reduce the rate of the perturbation to an
exponentially small quantity such that its square is much smaller than
the result. Then the Poincar\'e-Melnikov method is applicable.

It is natural to measure the rate of the separatrix splitting by
the area ${\cal A}_{\mbox{lobe}}$ of a lobe domain, bounded by
segments $I^{s,u}$ of the stable $(s)$ and unstable $(u)$ separatrix
such that $I^s$ and $I^u$ have the same boundary points, these
points are homoclinic, and there are no other common points of
$I^s$ and $I^u$.

Consider the system with Hamiltonian
\begin{equation}
\label{assum:small_period}
   \widehat H(\widehat x,\widehat y,t)
 = \varepsilon ( \widehat y^2 /2 + (1 + 2B\cos t) \cos{\widehat x} ),
\end{equation}
where $\widehat x,\widehat y$ are canonically conjugated variables
and $\varepsilon>0$ is small.
Theorem~1 implies that there exists a symplectic change
of coordinates $\widehat x, \widehat y \mapsto x, y$ which is

i) close to the identity,

ii) $2\pi$-periodic in time,

iii) real-analytic in a complex neighborhood of the separatrices
$\Gamma^\pm$ of the system with Hamiltonian $H_0 = \hat y^2 / 2 + \cos \hat x$

iv) such that the new Hamiltonian function takes the form:
$$
  H (x, y, t) = \varepsilon \big( H_0(x, y)
      + \varepsilon H_1 (x, y, \varepsilon)
      + \exp( - c / \varepsilon ) H_2 (x, y, t, \varepsilon)
                            \big) .
$$
Here $H_0 = y^2 / 2  +  \cos x$, the constant
$c \in [0,\pi/2)$ is arbitrary, the functions
$H_1, H_2$ are real-analytic in $x,y$ in the vicinity of the unperturbed
separatrices of the hyperbolic fixed point
$(\widehat x,\widehat y) = (0,0)$,
smooth in $\varepsilon > 0$ and the function $H_2$ is analytic
and $2\pi$-periodic in $t$.

The ordinary Poincar\'e-Melnikov theory applied to this system for positive
$\tilde c$ gives a correct asymptotics of the separatrix splitting,
\cite{tre_RJMP}. The following estimate holds:
$$
   {\cal A}_{\mbox{lobe}}
 = \frac{8\pi}{\varepsilon} \exp\Big(-\frac{\pi}{2\varepsilon}\Big)
   \Big( B f(B^2) + O\big(\frac{\varepsilon}{\log\varepsilon}\big) \Big),
$$
where $f$ is an entire real-analytic function, $f(0) = 2$.
In \cite{tre_RJMP, tre_cha} numerical values of several
Tailor coefficients of $f$ are presented.

Analogous results for the Standard Chirikov Map and for some Hamiltonian
systems with 2 degrees of freedom can be found in
\cite{tre_cha,tre_NATO,Nov}.

\subsection{Inclusion of a map into a flow}
\vskip-5mm \hspace{5mm}

In this section we consider the following problem:
to present a given self-map of a manifold $M$ as the time-one
map (the Poincar\'e map) in some ODE system, generated by a
periodic in time vector field.

The problem can be formulated for various classes of maps
and the corresponding vector fields.
For example, it is possible to consider generic maps and vector
fields, reversible ones with respect to some involution,
Hamiltonian, preserving a volume, etc.
Here we discuss the analytic set up i.e., assuming that
the map is real-analytic, we look for its inclusion into
a flow generated by a real-analytic vector field.
The problem in $C^\infty$ category is much simpler.

The following construction is well-known. Given a diffeomorphism $T$
of a manifold $M$ onto itself consider the direct product
$M\times [0,1]$ with the vector field $\partial/\partial t$, where
$t$ is the coordinate on $[0,1]$. The map $T$ generates the identification
$$
  M\times\{0\} \sim M\times\{1\},\quad
  (z,0) \sim (T(z),1).
$$
This identification converts $M\times [0,1]$ into a manifold
$\cal M$.
Let $\pi : M\times [0,1] \to {\cal M}$ be the natural projection.
The smooth vector field $\partial/\partial t$ generates on the surface
$\pi(M\times \{0\})\subset {\cal M}$
the Poincar\'e map coinciding with $T$.

This construction does not solve the problem we deal with because in
general it is not clear if $\cal M$ is real-analytically diffeomorphic
to $M\times{\bf T}^1$. Nevertheless, sometimes this can be proven.

The problem is solved in the symplectic set up for maps which are close to
integrable, \cite{dou2, kuk, kuk_poe} and for generic maps \cite{trif}.

Note that all the known proofs in the analytic set up use essentially
the Grauert theorem on the inclusion of an analytic manifold into
Euclidean space (or modified versions of this theorem).
Our result is based on the method of continuous averaging.

First, note that any map which is not isotopic to the identity\footnote
{Two smooth maps $T_j:M'\to M''$, $j = 0,1$
 ($M'$ and $M''$ are manifolds)
 are called isotopic if there exists a family of maps
 $\widehat T_s:M'\to M''$ of the same smoothness class
 continuous in the parameter $s\in[0,1]$, such that
 $\widehat T_0 = T_0$ and $\widehat T_1 = T_1$. In other words,
 if $T_0$ can be continuously deformed into $T_1$.}
obviously can not be included into a flow.
Let $M$ be an $m$-dimensional compact real-analytic manifold.
Let $\chi$ be a closed subalgebra in the Lie algebra
$(\cal L,\,[\,,])$ of all analytic vector fields on $M$.

We denote by $X$ the subset of all analytic diffeomorphisms of $M$ obtained
as a result of the time-one shift along solutions of a system
\begin{equation}
\label{shift}
        \dot z = u(z,t),\qquad  u(\cdot,t)\in\chi,
        \quad t\in[0,1],\quad z\in M
\end{equation}
(it is not assumed that $u(z,0) = u(z,1)$).
We assume that the vector field $u$ is $C^2$-smooth with respect
to time. This smoothness condition is technical. For example, it can
be replaced by continuity in $t$ in Hamiltonian case and
in the general one.

Obviously, all diffeomorphisms from $X$ are isotopic to the
identity inside $X$.

{\bf Theorem 3}
{\it For any map $T\in X$ there exists a vector field
$$
        U = U(z,t),\qquad  U(\cdot,t)\in\chi,\quad
        t\in {\bf R}, \quad  z\in M
$$
which is analytic in $z$ and $t$, $2\pi$-periodic in $t$, and such that
the time-$2\pi$ shift along its trajectories coincides with $T$.
}

As a corollary we obtain a possibility of the inclusion of analytic
maps into analytic flows in general, symplectic, and volume-preserving
cases.

The vector field $U$ is obviously not unique.

It can be also proved that if $T$
is reversible with respect to some involution $I:M\to M$,
$I^2 = \mbox{id}_M$ (i.e., $T\circ I = I\circ T^{-1}$), the
corresponding vector field $U$, can be also regarded $I$-reversible:
$U(z,t) = -dI\, U(Iz, -t)$.

Suppose that the map $T$ is close to $T_0$
($\mbox{\rm dist}(T,T_0) =\varepsilon$ in a complex neighborhood of $M$),\footnote
 {The distance can be defined for example as follows:
  $$
    \mbox{dist}(T,T_0) = \sup_{z\in M'} \rho(T(z), T_0(z)),
  $$
  where $M'$ is a complex neighborhood of $M$, and $\rho$ is some metric.
  The choice of the neighborhood and of the metric plays no role due to
  compactness of the manifold $M$.}
where $T_0$ is included into the flow generated by a periodic
analytic vector field $U_0$. Then the vector field $U$ can be chosen
close to $U_0$ ($|U-U_0| = O(\varepsilon)$
in a smaller complex neighborhood of $M$).
In particular, in the symplectic case if $T$ is close to an integrable map,
a Hamiltonian system associated with $T$ also can be chosen close to an
integrable one and the orders of closeness are the same.

Continuous averaging in the proof is used to smooth out the dependence
of the original vector field (\ref{shift}) on time, \cite{Pro_Tre}.

{\bf Acknowledgements.}
The work was partially supported by Russian Foundation of Basic Research
grants 02-01-00400 and 00-15-99269, and by INTAS grant 00-221.

\label{lastpage}

\end{document}